\documentclass[12pt]{article}
\usepackage{amsmath,amsxtra,amssymb,latexsym, amscd}
\usepackage[left=3cm,right=2.5cm,top=2.5cm,bottom=2.5cm]{geometry}

\newtheorem{theo}{Theorem}

\newcommand{\seq}[1]{\left<#1\right>}
\newcommand{\abs}[1]{\left\vert#1\right\vert}

\begin{document}
\newcommand{\lon}{\longrightarrow}
\newcommand{\lom}{\longmapsto}

\title{\bf  Analytic Extension of a maximal surface in $\Bbb L^3$ along its boudary\thanks{This work was supported by the National Basis Reseach Program 101706, Vietnam.}}          

\author{\bf Doan The Hieu and Nguyen Van Hanh\\
\sl  Dept. of Math., College of Education, Hue University\\
34 Le Loi, Hue, Vietnam\\
\sl deltic@dng.vnn.vn}        
\maketitle
\begin{abstract}
 We prove that a maximal surface in Lorentz-Minkowski space $\Bbb L^3$ can be extended analytically along its boundary if the boundary lies in a plane meeting the surface at a constant angle.
\end{abstract}

{\bf Subjclass:} {Primary 53C50; Secondary 58D10, 53C42}


\section{Introduction}      
A maximal surface in Lorentz-Minkowski space $\Bbb L^3$ is a spacelike surface with zero mean curvature. Maximal surfaces share many interesting properties with their counterparts, minimal surfaces, in $\Bbb R^3.$ For example, they are critical points (the maxima) of area variations and also admit Enneper-Weierstrass representations. It is well known that a minimal surface in $\Bbb R^3$ can be extended (symmetrically) along its boundary if the boundary lies in a plane meeting the minimal surface orthogonally. This fact also holds for maximal surfaces in $\Bbb L^3$ (see \cite{mi2}), where the plane is assumed to be timelike since spacelike and lightlike planes can not meet a maximal surface orthogonally, except at singular points, see the Remark in section 3.
\par
In 1996, J. Choe (\cite {choe}) proved that a minimal surface in $\Bbb R^3$ can be extended analytically along its boundary if the boundary lies in a plane meeting the minimal surface at a constant angle. The main idea is based on Enneper-Weierstrass representation of a minimal surface in terms of a holomorphic function $f$ and a meromorphic function $g.$ The meromorphic function $g$ can be viewed as the Gauss map of the minimal surface. Since the plane meets the minimal surfaces at a constant angle, the image of the boundary under the Gauss map $g$ lies in a circle and hence we can apply Schwartz reflection principle to extend both $f$ and $g$ along the boundary. 

In this paper, we show that the above idea can be apply for the case of maximal surfaces in $\Bbb L^3.$ The complication in this stuation is a  plane can be spacelike, timelike or lightlike. 


\section{Preliminaries}
 The Minkowski 3-space $\Bbb L^{3}$ is the 3-dimensional vector space $\Bbb R^{3}=\{(x_1, x_2, x_3,): x_i\in \Bbb R, i=1,2,3\}$ endowed with the indefinite $(2,1)$-metric 
$$\langle x, y\rangle=x_1y_1+x_2y_2-x_3y_3,$$
where $x=( x_1, x_2, x_3), y=( y_1, y_2, y_3)\in \Bbb L^{3}.$

 We say that a nonzero vector $x\in \Bbb L^{3}$ is \emph{spacelike}, \emph{lightlike} or \emph{timelike} if  $\langle x, x\rangle>0$, $\langle x, x\rangle=0$ or $\langle x, x\rangle<0$, respectively. The vector zero is always considered as a spacelike one.

The norm of a vector $x\in \Bbb L^{3}$, denoted by $\|x\|$, is defined by $\sqrt{|\langle x,x\rangle|}.$ The definition of the cross-product of two vectors $a=(a_1,a_2,a_3);\ b=(b_1,b_2,b_3),$ denoted by $a\wedge b$  is given as follow
$$a\wedge b=(a_2b_3-a_3b_2, a_3b_1-a_1b_3, a_2b_1-a_1b_2).$$

For a nonzero vector $n\in \Bbb L^{3}$, a plane with (pseudo) normal $n$ is the set
$$P(n,c)=\{x\in\Bbb L^{3} : \langle x,n\rangle=c, c\in\Bbb R\}.$$

The plane $P(n,c)$ is called \emph{spacelike}, \emph{lightlike} or \emph{timelike} if $n$ is timelike, lightlike or spacelike, respectively.

It is easy to see that,  $P(n,c)$ is spacelike if any vector $x\in P(n,c)$ is spacelike; $P(n,c)$ is lightlike if  $P(n,0)$ is tangent to the lightcone; $P(n,c)$ is timelike if it contains timelike vectors.

The set
$$\Bbb H^2=\{x\in \Bbb L^3 : \langle x,x\rangle =-1\}$$
is called the Hyperbolic. It has two connected components $\Bbb H^2_+=\{x\in \Bbb H^2 : x_3 \ge1\}$ and  $\Bbb H^2_-=\{x\in \Bbb H^2 : x^3\le-1\}.$
For studying spacelike surfaces in Lorentz-Minkowski spaces, $\Bbb H^2,\ \Bbb H^2_-, \ \Bbb H^2_+$ play the same roles as the unit sphere $\{|x|^2=1\}$ in Euclidean spaces.

Let $X:M \lon\Bbb L^{3}$ be an immersion of a 2-dimensional connected manifold.  X (or $X(M)$) is called spacelike if the induced metric on $M$ via $X$ is a Riemannian metric. That means the tangent plane $T_pM\subset T_p\Bbb L^{3}$
 is spacelike, for every $p\in M.$ In this case, the manifold $M$ is orientable.
Now suppose that $X:M \lon\Bbb L^{3}$ is a spacelike immersion and $(u,v)$ is a local coordinate system. The (local) unit normal vector field is defined as follow
$$N(u,v)=\frac{X_u\wedge X_v}{\|X_u\wedge X_v\|}.$$ 
Because $M$ is spacelike, $N(u,v)$ is always timelike. 

Since $M$ is connected, we can define the unit normal timelike vector field $N$ on $M$ globally and the image of $N$ lies in one of  components of $\Bbb H^2.$  Because of that we   can consider $N$  as a map $N:M\lon \Bbb H^2_+.$ The map $N$ is  called the Gauss map of the immersion. The shape operator is the map $A:=-dN$ defined for all vector fields on manifold $M$ and the mean curvature $H$ is a half of the trace of $A$
$$ H:=\frac 12 tr(A).$$

A spacelike immersion $X:M \lon\Bbb L^{3}$  is said to be a maximal immersion if its mean curvature (at every point) is equal to zero, that is $H=0.$

In 1983, Kobayashi (\cite{ko1}) showed  Enneper-Weierstrass representations for maximal immersions in $\Bbb L^3.$ Such representations for maximal immersions   are quite similar as that for minimal immersions in Euclidean space $\Bbb R^3.$ It is clear that, we can defined a local isothermal coordinate systems, whose changes of coordinates preserve the orientation, for  maximal immersions. The existence of such coordinate systems is proved quite similar  as that for minimal immersions. Thus, since every spacelike immersion is orientable, $M$ admit a structure of a Riemann surface.  

Now suppose $X=(x_1,x_2,x_3)$ and $z=u+iv$ be the local complex parameter on $M.$
 We set 
$$\phi_k:=\frac 12\left(\frac {\partial x_k}{\partial u}-i\frac {\partial x_k}{\partial v}\right),\ \ \ k=1,2,3.$$
Since $M$ is maximal, $x_k,\ k=1,2,3$ are harmonic and hence $\phi_k,\ k=1,2,3$ are holomorphic. 
Direct computation shows that 
\begin{equation}\label{21}      \phi_1^2+\phi_2^2-\phi_3^2=0,  \end{equation}
and
\begin{equation}\label{22} |\phi_1|^2+|\phi_2|^2-|\phi_3|^2>0. \end{equation}
We see that  $ds^2=|\phi_1|^2+|\phi_2|^2-|\phi_3|^3>0$ is the Riemannian metric on $M$ induced by the immersion $X$ and  $\phi_k,\ k=1,2,3$ have no real periods and hence the immersion $X$ can be represented as
\begin{equation}\label{23} X(z)={\text{Re}}\int (\phi_1,\phi_2,\phi_3)dz, \end{equation}
where the integral is taken on  an arbitrary path from a fixed point to $z.$ 

Conversely,  if $\phi_1, \phi_2,\phi_3$ are holomorphic functions on $M$ that have no real periods and satisfy (\ref{21}) and (\ref{22}), then (\ref{23}) determines a maximal surface. 

If $\phi_1-i\phi_2=0,$ then $\phi_3=0.$ In this case, $M$ is a plane.
Now suppose that $\phi_1-i\phi_2\ne 0,$ we set
$$f=\phi_1-i\phi_2,$$
$$g=\dfrac{\phi_3}{\phi_1-i\phi_2}.$$ 
We have 
\begin{equation}\label{24}
\begin{cases}            \phi_1=\frac12f(1+g^2)\\
\phi_2=\frac{i}{2}f(1-g^2)\\
\phi_3=fg \end{cases}, \end{equation} 
and then (\ref{23}) can be writen as follow
 \begin{equation}\label{25}
X(z)=\mathrm{Re}\left(\frac12\int_{z_0}^zf(1+g^2)d\omega,\dfrac{i}{2}\int_{z_0}^zf(1-g^2)d\omega,\int_{z_0}^zfgd\omega\right).
\end{equation}
From (\ref{24}),  we have $\phi_1+i\phi_2=fg^2.$ Thus, we can conclude that the poles of $g$ coincide with the zeroes of $f$ in such a way that a pole of order $2m$ of $g$ corresponds to a zero of order $m$ of $f.$  Conversely,  if such $g$ and $f$ is given then (\ref{25}) determines a maximal immersion.

Since  $(X_u-iX_v)=2(\phi_1,\phi_2,\phi_3),$ we have 
$$\begin{aligned}                
X_u\wedge X_v
&=4\mathrm{Im}(\phi_2\overline{\phi}_3,\phi_3\overline{\phi}_1,\phi_2\overline{\phi}_1)\\
&=\abs{f}^2(1-\abs{g}^2)\left(2\mathrm{Re}(g),2\mathrm{Im}(g),
1+\abs{g}^2\right).  \end{aligned}$$ 
Thus, the Gauss map $N$ can be expressed as follow
$$
N=\left(\dfrac{2\mathrm{Re}(g)}{1-\abs{g}^2},
\dfrac{2\mathrm{Im}(g)}{\abs{1-g}^2},\dfrac{1+\abs{g}^2}{1-\abs{g}^2}\right).
$$
Since  $N(z)\in \Bbb H^2_+,$ we conclude that $|g|<1.$

It is clear that $z\lom \left(\dfrac{2\mathrm{Re}(z)}{1-\abs{z}^2},
\dfrac{2\mathrm{Im}(z)}{1-\abs{z}^2},\dfrac{1+\abs{z}^2}{1-\abs{z}^2}\right)$ is a conformal isomorphism $\pi$ between $ D=\{z\in\Bbb C: |z|<1\}$ and $\Bbb H^2_+.$  The map $\pi^{-1}$ is the stereographic projection from the point $(0,0,-1).$ The formular of $\pi^{-1}$ is
$$\pi^{-1}(x_1,x_2,x_3)=\frac{x_2+ix_2}{1+x_3}.$$ We can view $g$ as a map from $M$ into $ D$ and $\pi^{-1}\circ N=g.$  Because of that we also call $g$ the Gauss map of $M.$ 
\section{Extension of a maximal surface}
 Let $\Omega$ be a domain in $\Bbb R^2.$ We will call a maximal immersion $x:\Omega\lon \Bbb L^3$ a maximal surface and always assume that the parameters $u,v$ on $\Omega$ are isothermal and set $z=u+iv.$.

It is well known that, every maximal immersion can be locally writen as a maximal surface and by Uniformization theorem a simply connected maximal immersion can be expressed as a maximal surface globally. 

Denote $D=\{u^2+v^2<1\},$ $D_+=\{u^2+v^2<1;\ v>0\},$ $D_-=\{u^2+v^2<1;\ v<0\}$ and  $D_0=D\cap \{v=0\};$ we have the main theorem of this paper.
\begin{theo}\label{th1}
Let $X_+:D_+\lon\Bbb L^3$ be a maximal surface with isothermal parameters $u,v$ and $\Pi$ be a plane. Suppose that $\gamma$ is an analytic
curve in $\Pi$, $X_+(u,v)$ tend to
$\gamma(u)$ whenever $v\rightarrow 0,$   and
$$\lim_{v\rightarrow 0}{\seq{N(z), n}}= c\ne 0,$$ 
where $N$ is the Gauss map of $X_+$ and $n$ is the unit normal vector of $\Pi.$ Then $X_+$ can be analytically extended along $\gamma$ to a maximal surfaces $X:D\lon \Bbb L^3$ such that $X|_{D_+}=X_+$ and $X(D_0)=\gamma.$
\end{theo}
{\bf Proof.}
The main idea for the proof is  showing that both $g$ and $f$ can be   extended analytically  on $D$ and hence by (\ref{25}) we get the extended maximal surface. We will consider three cases: $\Pi$ is spacelike, $\Pi$ is timelike and $\Pi$  is lightlike.
In each case, we will use the following fact: if $g$ can be continuously extended to $D_+\cup D_0$ and $g(D_0)$ lies in a circle, then after using a M\"{o}bius transformation that maps $g(D_0)$ to the real axis, Schwartz reflection principle can be applied to extend $g$ on $D.$ 
\par
\begin{enumerate} 
\item {\bf $\Pi$ is spacelike.}
By using a suitable Lorentzian transformation,  we can assume that $\Pi$ is the plane $x_3=0.$ We choose $n$ is timelike vector $(0,0,1),$  then

$$\lim_{v\to0}{\seq{\xi(z),n}}=\lim_{v\to0}\dfrac{1+\abs{g}^2}{1-\abs{g}^2}=-c.$$ 
Set $c= \cosh\theta,$ we conclude that

$$
\lim_{v\to0}\abs{g(z)}=\coth\dfrac{\theta}{2}.
$$
 The meromorphic function $g$ then can be continuously extended on $D_+\cup D_0$ such that
$$\abs{g(z)}=\coth\dfrac{\theta}{2},\ \ \forall z\in  D_0.$$
That mean $g$ maps $D_0$ into the circle with the center  $O$ and radius 
$r=\coth\dfrac{\theta}{2},$ and therefore, $g$ can be extended analytically on $D.$ The extension of $g$ also denote by $g$ and is expressed as follow:
$$
g(z)=\coth^2(\dfrac\theta2)(\overline{g(\overline{z})})^{-1},z\in D_-.
$$

Next, we extend $f$ on $D.$ First we observe that $x_3$ extends to a harmonic function, also denoted by $x_3$, on $D$ by setting 
$$
x_3(z)=-x_3(\overline{z}),\ \ z\in D_-.
$$
Then $\phi_3$ can be extended to a holomorphic function, also denoted by $\phi_3$, on $D$ by setting 
$$
\phi_3(z)=-\overline{\phi_3(\overline{z})},z\in\Omega_-.
$$
Finally, $f$ is extended analytically on $D$ by setting
$$
f(z)= \dfrac{-\overline{f(\overline{z})g^2(\overline{z})}}{\coth^2(\dfrac\theta2)}=
\dfrac{-\overline{\phi_3(\overline{z}){g(\overline{z})}}}{\coth^2(\dfrac\theta2)}
=\dfrac{-\overline{\phi_3(\overline{z})}}
{\coth^2(\dfrac\theta2)(\overline{g(\overline{z})})^{-1}}=\dfrac{\phi_3(z)}{g(z)} 
,\ \forall z\in D_-.
$$
\par
\item{\bf $\Pi$ is timelike.}

We can assume that $\Pi$ is the plane $x_2=0.$ Set $c=\dfrac{1}{\lambda},\ \lambda\ne 0$ and choose $n=(0,1,0).$ 
The assumption
$\lim\limits_{v\to0}{\seq{N(z),n}}=\dfrac{1}{\lambda},$ implies that
$$\lim_{v\to0}\dfrac{2\mathrm{Im}(g)}{1-\abs{g}^2}=\dfrac{1}{\lambda}.$$

Thus, $g$ is extended continuously on $D_+\cup D_0$ such that the following is satisfied
\begin{equation}  \label{g2}  
\dfrac{2\mathrm{Im}(g)}{1-\abs{g}^2}=\dfrac{1}{\lambda}.
  \end{equation} 
Equation (\ref{g2}) gives 
 $$[\mathrm{Re}(g)]^2+[\mathrm{Im}(g)+\lambda]^2=1+\lambda^2.$$
Therefore, $g$ maps $D_0$ into the circle with the center at $(0,-\lambda)$ and radius $r=\sqrt{1+\lambda^2}$ and then the meromorphic function  $g$ is extended as follow
$$
g(z)=-i\lambda +(1+\lambda^2)\left(\overline{g(\overline{z})}-i\lambda
\right)^{-1},\forall z\in D_-.$$

Because $X_+=(x_1,x_2,x_3)$ is maximal and $u,v$ are isothermal parameters, $x_2$ is harmonic on $D_+.$ The assumption
$$\lim_{v\to0}X_+(z)=\gamma(u)\in \Pi,$$
implies that $x_2$ can be continuously extended on $D_+\cup D_0$ by setting
$$x_2(z)=0,\forall z\in D_0.$$ 

Schwartz reflection principle says that $x_2$ can  be extended on $D$ as follow
$$
x_2(z)=-x_2(\overline{z}),\forall z\in D_-.
$$
Therefore, $\phi_2$ is extended on $D$ by setting
$$
\phi_2(z)=-\overline{\phi_2(\overline{z})}, z\in \Omega_-.
$$
Since  $g^2(z)\ne 1,$ the holomorphic $f$ is extended analytically on $D$ by setting
$$
f(z)=\dfrac{2\phi_2(z)}{i(g^2(z)-1)},z\in D.
$$
\par
\item {\bf $\Pi$ is lightlike.}

Assume that the equation of 
 $\Pi$ is $x_1-x_3=0.$ 
We set $c=1+\lambda$ and choose $n=(1,0,1).$

If $\lambda=0,$ then by the assumption
$$\lim\limits_{v\to 0}{\seq{N(z),n}}=1,$$ we have
$$\lim_{v\to 0}\left[ \dfrac{2\mathrm{Re}(g)}{1-\abs{g}^2}-\dfrac{1+\abs{g}^2
}{1-\abs{g}^2} \right]=1$$
  or equivalently,
$\mathrm{Re}(g)$ tends to $1$ whenever $v$ tends to $0.$
The meromorphic function $g$ can be  extended continuously on $D\cup D_0$ such that $\mathrm{Re} g(z)=1,\ \ \forall z\in D_0.$
That mean $g(D_0)\subset \{\mathrm{Re}(z)=1\}.$ 
In this case $g$ is extended analytically on $D$ by setting
$$
g(z)=2-\overline{g(\overline{z})},\forall z\in D_-.
$$
If $\lambda\ne 0,$ by the assumption  
$$\lim_{v\to0}{\seq{N(z),n}}=
\lim_{v\to0}\left[ \dfrac{2\mathrm{Re}(g)}{1-\abs{g}^2}-\dfrac{
1+\abs{g}^2}{1-\abs{g}^2} \right]=1+{\lambda},
$$
we conclude that  $g$ can be continuously extended on $D\cup D_0$ in such away that
$$
\left( \mathrm{Re}(g)+\frac 1\lambda \right)^2+\left(\mathrm{Im}(g) \right
)^2=(1+\frac 1\lambda)^2,\forall z\in D_0.
$$
That means the image of  $D_0$ under  $g$ lies in the circle with center
$(-\frac 1\lambda,0)$ and radius $r=\abs{1+\frac 1\lambda}.$

Then  $g$ is extended analytically on $D$ by setting
$$
g(z)=-\frac 1\lambda+(1+\frac 1\lambda)^2\left( \overline{g(\overline{z})}+\frac 1\lambda\right
)^{-1},\forall z\in D_-.
$$
In order to extend $f$ we first observe that
 $\psi=x_1-x_3$ is a harmonic function on $D_+$  and by the assumption
$$\lim_{v\to0}X_+(z)\to\gamma(u)\in \Pi,$$ 
it can be extended to a continuous function  on  $D_+\cup D_0$ by setting $\psi(z)=0,\forall
z\in\Omega_0.$ Then by Schwartz reflection principle for harmonic function, $\psi$ can be extended to a harmonic function  on  $D$ by setting
$$
\psi(z)=-\psi(\overline{z}),\forall z\in D_-.
$$

Let $\psi^*$ be the harmonic conjugation of $\psi$ then 
$\dfrac{d(\psi+i\psi^*)}{dz}$ is a holomorphic function on $D.$ It is clear that $\dfrac{d(\psi+i\psi^*)}{dz}\Big|_{D_+}=\phi_1-\phi_3.$ So $\dfrac{d(\psi+i\psi^*)}{dz}$ is the extension of $\phi_1-\phi_3 $ and we can write $\phi_1-\phi_3 $
instead of $\dfrac{d(\psi+i\psi^*)}{dz}.$
Then, the analytic extension of $f$ can be writen as follow
$$f(z)=\dfrac{2(\phi_1-\phi_3)}{(1-g(z))^2},
\forall z\in D_-.$$ 
    \end{enumerate} 
{\bf Remark.} 
\begin{enumerate} 
\item It is clear that  if $z=u+iv$ is a pole of order $m$ of  $g$  then $z=u+iv$ is a  zero of  order $2m$ of $f.$
\item If $\langle N(z), n\rangle =0 $ along $S\cap\Pi,$ we then say that the plane $\Pi$ meets the maximal surface $S$ orthogonally. Suppose $\Pi$ is spacelike, then we can conclude that $1+|g|^2=0,$ a contradiction. Thus, a spacelike plane can not meet a maximal surface orthogonally. 
If $\Pi$ is lightlike, we can suppose the equation of $\Pi$ is $x_1-x_3=0.$ Then we can conclude that $g=-1$ along  $S\cap\Pi.$  Therefore, $X_u\wedge X_v=0.$ Thus a lightlike plane can not meet a maximal surface orthogonally, except at singular points.. 
\item We can see the extension clearly on Lorentzian Catenoid. Let 
$$X(u,v)=(\sinh u\cos v, \sinh u\sin v, u);\ \ (u,v)\in U=\{(u,v)\in \Bbb R_+\times(-\pi,\pi) $$
be the   Lorentzian Catenoid with only conelike sigularity at the origin. Let $\Pi_1$ be spacelike plane $x_3=a>0,$ $\Pi_2$ be spacelike plane $x_3=b>a $ and $\Pi_3$ be spacelike plane $x_3=2b-a. $  The  extension about $\Pi_2$ as in proof of Theorem \ref{th1} maps $X(U)\cap\Pi_1$ to  $X(U)\cap\Pi_3$  and maps the component bounded by $\Pi_1$ and $\Pi_2$ to the component bounded by $\Pi_2$ and $\Pi_3.$ 
\item ({\bf Extension about a conelike singularity})
   Nevertheless, there are important differences between maximal surfaces and minimal surfaces. The fact that the only complete maximal surfaces in $\Bbb L^3$ are spacelike planes is an example in global theory. On the other hand, maximal surfaces may have isolated singularities that never happen for minimal surfaces (see \cite{ko2}).

Let $S$ be a maximal surface and $p\in S$ is a conelike singularity. For more detail about conelike singularities, we refer the readers to \cite{ko2}. Now let $X:\overline D\lon \Bbb L^3$ be a neighbourhood of a conelike singularity where $X(0,0)$ is the conelike singularity and suppose that $X(\partial D)$ meets spacelike plane $x_3=a$ at a constant angle. In this situation, the image of $g$ is an annulus bounded by circles $\{|z|=1\}$ and $\{|z|=r<1\}$ and hence conformally identified  with $D-\{(0,0)\}.$ Obviously, we can extend both $\phi_3$ and $g$ analitycally to the whole $\Bbb C$ by using the inversion about circle $\{|z|=1\}.$  The resulting is a complete maximal surface with one conelike singular point and one end and therefore is an embedding entire graph (see \cite{fe}, Proposition 2.1). It must be the Lorentzian Catenoid by a result of Ecker (see \cite{ec}).

\item  ({\bf Extension about an end}) The same argument as in item 4 also holds for $X:\overline D-\{(0,0)\}\lon \Bbb L^3$ being a neighbourhood of an end of a maximal surface, and $X(\partial D)$ meets spacelike plane $x_3=a$ at a constant angle. In this case, the image of gauss map $g$ is the disk $\{0<|x|<r;\ r<1\}$ and also can be extended analitycally to $D-\{(0,0)\}.$

 \end{enumerate}

\end{document}